\documentclass[12pt]{article}

\usepackage{amsmath,amsfonts,amssymb,verbatim}
\newcommand{\supp}{\operatorname*{supp}}

\newcommand{\be}{\begin{equation}}
\newcommand{\ee}{\end{equation}}
\newcommand{\la}{\label}
\newcommand{\ba}{\begin{array}{c}}
\newcommand{\ea}{\end{array}}

\newcommand{\nax}{\nabla_x}

\newtheorem{thm}{Theorem}[section]

\newtheorem{test}{Definition}[section]

\newtheorem{sssvort}[test]{Definition}
\newtheorem{testII}[test]{Definition}
\newtheorem{boundsupp}[test]{Proposition}
\newtheorem{Eulerdef}[test]{Definition}
\newtheorem{Prokhorov}[test]{Theorem}
\newtheorem{ilimit}[test]{Theorem}
\newtheorem{corProk}[test]{Theorem}
\newtheorem{Glimit}{Definition}[section]

\newtheorem{lemma}{Lemma}[section]

\title{Inviscid limit for damped and driven incompressible Navier-Stokes 
equations in ${{\mathbb R}^2}$}
\author{P. Constantin\\Department of Mathematics, The University of Chicago\\ Chicago IL 60637 USA,\\F. Ramos \\Instituto de Matem\'atica, Universidade Federal \\do Rio de Janeiro\\Rio de Janeiro RJ 21945-970 Brazil.}

\begin{document}

\maketitle

\noindent{\bf Abstract.} We consider the zero viscosity limit of long time
averages of solutions of damped and driven Navier-Stokes equations in ${\mathbb R}^2$. We prove that the rate of dissipation of enstrophy vanishes. Stationary statistical solutions of the damped and driven Navier-Stokes equations converge to renormalized stationary statistical solutions of the damped and driven Euler equations. These solutions obey the enstrophy balance.


\noindent{\bf key words} Inviscid limit, statistical solutions, anomalous dissipation.

\section{Introduction}
The vanishing viscosity limit of solutions of 
 Navier-Stokes equations is a subject that has been
extensively studied. Boundary layers, which present the 
most important physical aspects of the problem, are difficult to study 
and their mathematical understanding is rather limited.
More progress has been made in the study of the limit when boundaries are absent (flow in ${\mathbb R}^n$ or ${\mathbb T}^n$). Even in this restricted situation, there are two distinct concepts of vanishing viscosity limit. The finite time, 
 zero viscosity limit is the limit $\lim_{\nu\to 0}S^{\nu}(t)(\omega_0)$ of
 solutions $S^{\nu}(t)(\omega_0)$ of the Navier-Stokes equations with a fixed initial datum $\omega_0$ and with time $t$ in some finite interval $[0,T]$. 
By contrast, in the infinite time zero viscosity limit, 
long time averages of functionals of the
solutions
$$
\lim_{t\to\infty}\frac{1}{t}\int_0^t\Phi(S^{\nu}(s)\omega_0)ds = \int\Phi(\omega)d\mu^{\nu}(\omega)
$$
are considered first, at fixed $\nu$. These are represented by measures $\mu^{\nu}$
in function space. The long time, zero viscosity limit is then 
$\lim_{\nu\to 0}\mu^{\nu}$, 
$$
\lim_{\nu\to 0}\left (\lim_{t\to\infty}\frac{1}{t}\int_0^t\Phi(S^{\nu}(s)\omega_0)ds\right).
$$
The two kinds of limits are not the same. This is most clearly seen in the 
situation of two dimensional, unforced Navier-Stokes equations. In this case,  any smooth solution of the Euler equations is a finite time inviscid limit but the infinite time inviscid limit is unique: it is the function identically equal to zero.  This simple example points out the fact that the
infinite time zero viscosity limit is more selective. In less simple situations, when the Navier-Stokes equations are forced, the long time inviscid limit is
not well understood. 

The finite time zero viscosity limit is the limit that has been most studied.
For smooth solutions in ${\mathbb {R}}^3$, the zero viscosity limit is given by solutions of the Euler equations, for short time, in classical (\cite{swann}),
and Sobolev (\cite{kato}) spaces; the limit holds for as long as the Euler solution is smooth (\cite{c0}). The convergence occurs in the  Sobolev space
$H^s$ as long as the solution remains in the same space (\cite{masmoudi}). The
rates of convergence are optimal in the smooth regime, $O(\nu)$. In some nonsmooth regimes (smooth vortex patches), the finite time inviscid limit exists and optimal rates of convergence can be obtained (\cite{abidi}, \cite{masmoudi}) 
but the rates deteriorate when the smoothness of the initial data deteriorates
 -- for nonsmooth vortex patches ({\cite{cwu1}}).

One of the most fundamental questions concerning the inviscid limit is: what happens to ideally conserved quantities?  For instance, in three dimensions, the kinetic energy is conserved by smooth Euler flow, and dissipated by viscous Navier-Stokes flow. Does the rate of dissipation of kinetic energy vanish with viscosity, or 
is there a non-zero limit? This is the problem of anomalous dissipation. The  term was coined relatively recently by field theorists but the anomaly was
suggested by Onsager and Kolmogorov independently in the nineteen forties. 
The problem is open. 

In two dimensions there exist infinitely many integrals that are conserved by smooth Euler flows. One of them is the enstrophy
$$
\int_{{\mathbb R}^2}\left |\omega (x,t)\right|^2dx
$$
where $\omega$ is the vorticity of the flow. The existence 
of anomalous dissipation of enstrophy is postulated in Kraichnan's theory 
for two dimensional turbulence (\cite{frisch}). 
This was studied in the framework of finite time inviscid limits with rough initial data (\cite{eyink}, \cite{lml}). It was established that, if the initial vorticity belongs to $L^2({\mathbb R}^2)$
then rate of dissipation of enstrophy vanishes with viscosity, for finite time. The finite time inviscid limits are weak solutions of Euler equations.

In this paper we study the long time, zero viscosity limit for damped and driven two dimensional Navier-Stokes equations. The damped and driven two dimensional equations arise in the Charney-Stommel model of the Gulf Stream (\cite{bct}). The fact there is no anomalous dissipation of energy in damped and driven Navier-Stokes equations was suggested by D. Bernard (\cite{bernard}). 

The paper is organized as follows. In the second section we describe the equations and a few of the properties of individual solutions of the viscous equations $S^{NS,\gamma}(t)(\omega_0)$. One of the facts that plays a significant role in the paper is that the positive semiorbit $O^{+}(t_0) = \{S^{NS,\gamma}(t)(\omega_0)\, \left |\,\,\, t\ge t_0 > 0\right.\}$ is relatively compact in $L^2({\mathbb R}^2)$ and included in a bounded set in $L^{1}({\mathbb R}^2)\cap L^{\infty}({\mathbb R}^2)$ that does not depend on the viscosity. The uniform bound uses essentially the fact that the damping factor $\gamma>0$ is bounded away from zero independently of the vanishing viscosity. In order to prove compactness,
because we work in the whole space, we need to prove also that the solution does not travel. Our results apply to the spatially 
periodic boundary conditions as well. The absence of anomalous dissipation of energy follows immediately from the bounds in the second section. 

The third section is devoted to the study of the vanishing viscosity 
limit of sequences of time independent individual solutions. The sequences have enough compactness to pass to convergent subsequences. The resulting solution is a weak solution of the damped  and driven Euler equations. The existence of weak solutions of such equations in the case of the Charney-Stommel model was first obtained in (\cite{bct}).  The weak solution of the damped and driven Euler equation is a renormalized solution in the sense of (\cite{dpl}). This implies that the weak solution obeys an enstrophy balance and that is used to show that there is no anomalous dissipation. 

The fourth section introduces the notion of stationary statistical solution of the damped and driven Navier-Stokes equations in the spirit of ({\cite{foias72}}, \cite{foias73}). In the case of finite dimensional dynamical systems 
$
\frac{d\omega}{dt} = N(\omega),
$
invariant measures $\mu$ obey
$\int \nabla_{\omega}\Psi(\omega)N(\omega)d\mu(\omega) = 0$ for any
test function $\Psi$. In infinite dimensions we need to restrict the
test  functions to a limited class of admissible functions. Among them are
generalizations of the characters $\exp{i\langle \omega, {\mathbf{w}}\rangle}$
with ${\mathbf {w}}$ a test function and an additional type of test function 
$\Psi_{\epsilon}(\omega)$ that uses $(\beta(\omega_{\epsilon}))_{\epsilon}$,
a mollification of a function of a mollification of $\omega$. Such technical 
precautions aside, the notion of  stationary statistical solution of the damped and driven Navier-Stokes equation is a natural extension of the notion  of  invariant measure for deterministic finite dimensional dynamical systems. We show that weak limits of stationary statistical solutions of the damped and driven Navier-Stokes equations are renormalized stationary statistical solutions of the damped and driven Euler equations, a concept that we introduce in the spirit of (\cite{dpl}). We also show that if the supports of the stationary statistical solutions of the damped and driven 
Navier-Stokes equations are included in sets that are bounded uniformly in in $L^p({\mathbb R}^2)\cap L^{\infty}({\mathbb R}^2)$ (with $p<2$ for technical reasons having to do with the slow decay at infinity of velocity in the Biot-Savart law) then the weak limits are renormalized stationary statistical solutions of the damped and driven Euler equations that obey the enstrophy balance.

In the fifth section we prove our main results. We construct stationary statistical solutions $\mu^{\nu}$ of the damped and driven Navier-Stokes equations by the Krylov-Bogoliubov procedure of taking long time averages. We show that
these solutions have good enough properties so that their weak limits are renormalized stationary statistical solutions $\mu^{0}$ of the damped and driven Euler equations that obey the enstrophy balance. We use this fact to prove that zero viscosity limit of the long time average enstrophy dissipation rate vanishes: 
$$
\lim_{\nu\to 0}\nu\left (\lim\sup_{t\to\infty}\frac{1}{t}\int_0^t \|\nabla \omega(s+t_0)\|_{L^2({\mathbb R}^2)}^2ds\right) = 0
$$
holds for all solutions $\omega (t) = S^{NS,\gamma}(t)(\omega_0)$, all $t_0>0$,
and all $\omega_0\in L^p({\mathbb R}^2)\cap L^{\infty}({\mathbb R}^2)$. We also prove that convergence in this class of statistical solutions is such that
$$
\lim_{\nu\to 0}\int_{L^2({\mathbb R}^2)}\|\omega\|^2_{L^2({\mathbb R}^2)}
d\mu^{\nu}(\omega) = \int_{L^2({\mathbb R}^2)}\|\omega\|^2_{L^2({\mathbb R}^2)}
d\mu^{0}(\omega). 
$$

\section{The setup}

We consider damped and  driven Navier-Stokes equations in ${\mathbb R}^2$
\be
\left \{
\ba
\partial_t u + u\cdot\nabla u -\nu\Delta u +  \gamma u + \nabla p = f,\\
\nabla\cdot u = 0
\ea
\right.
\la{ueq}
\ee
with $\gamma>0$ a fixed damping coefficient, $\nu>0$, $f$ time independent  with zero mean and $f\in W^{1,\infty}(\mathbb{R}^2)\cap H^1(\mathbb{R}^2)$. The initial velocity is divergence-free and belongs to $(L^2({\mathbb R}^2))^2$. We start by stating some of the properties of the individual solutions.
\begin{thm} \label{indi} Let $u_0$ be divergence-free, $u_0\in H^1({\mathbb R}^2)^2$. Then the solution of (\ref{ueq}) with initial datum $u_0$ exists for all time, is unique, smooth, and obeys the energy equality 
\be
\frac{d}{2dt}\int_{{\mathbb{R}^2}} |u|^2 dx + \gamma\int_{{\mathbb R}^2}|u|^2dx+\nu\int_{{\mathbb R}^2}|\nabla u|^2 dx
 = \int_{{\mathbb R}^2}f\cdot u dx.
\la{enba}
\ee
The kinetic energy is bounded uniformly in time, with bounds independent of
viscosity:
$$
\|u(\cdot,t)\|_{L^2({{\mathbb R}^2})} \le e^{-\gamma t}\left\{\|u(\cdot,0)\|_{L^2({{\mathbb R}^2})} - \frac{1}{\gamma}\|f\|_{L^2({{\mathbb R}^2})}\right\} + \frac{1}{\gamma}\|f\|_{L^2({{\mathbb R}^2})}.
$$
The vorticity $\omega$ (the curl of the incompressible two dimensional velocity)
\be
\omega = \partial_1 u_2 - \partial_2 u_1 = \nabla^{\perp}\cdot u
\la{curlu}
\ee
obeys
\be
\partial_t \omega + u \cdot\nabla \omega -\nu\Delta\omega  +\gamma\omega = g,
\la{omeq}
\ee
with   $g \in L^2({\mathbb R}^2)$, the  vorticity source,
 $g = \nabla^{\perp}\cdot f$. The map $t\mapsto \omega(t)$ is continuous $[0,\infty)\to L^2({\mathbb R}^2)$. If the initial vorticity is in $L^p({\mathbb R}^2)$, $p\ge 1$, and $g\in L^p({\mathbb{R}}^2)$, then the $p$-enstrophy is bounded uniformly in time, with bounds independent of viscosity
$$
\|\omega(\cdot,t)\|_{L^p({{\mathbb R}^2})} \le  e^{-\gamma t}\left\{\|\omega(\cdot,0)\|_{L^p({{\mathbb R}^2})} - \frac{1}{\gamma}\|g\|_{L^p({{\mathbb R}^2})}\right\} + \frac{1}{\gamma}\|g\|_{L^p({{\mathbb R}^2})}
$$
for $p\ge 1$. 
Moreover, the solution does not travel:
For every $\epsilon>0$ there exists $R>0$ such that, 
$$
\int_{|x|\ge R}\left | \omega(x,t)\right|^2 dx \le \epsilon
$$
holds for all $t\ge 0$.
\end{thm}
The proof of this theorem uses well-known methods, and will not be presented here. We only sketch the proof of the last statement. We take a smooth nonnegative function $\phi$ supported in $\{x\in {\mathbb {R}^2}\, |x|\ge \frac{1}{2}\}$ and identically equal to $1$ for $|x|\ge 1$, multiply the vorticity equation (\ref{omeq}) by $\phi\left({\frac{x}{R}}\right )\omega(x,t)$ and integrate in space. Denoting
$$
Y_R(t) = \int\phi\left (\frac{x}{R}\right )\left|\omega(x,t)\right|^2dx
$$
we obtain:
$$
\ba
\frac{d}{2dt} Y_R(t) + \gamma Y_R(t) \le \\
C\left\{\sqrt{Y_R(t)\int_{|x|\ge\frac{R}{2}}|g(x)|^2dx} + \frac{\nu}{R^2}\int|\omega(x,t)|^2dx  +
\frac{1}{R}\int |u(x,t)||\omega(x,t)|^2dx\right\}
\ea
$$
We deduce that
$$
\ba
\frac{d}{dt} Y_R(t) + \gamma Y_R(t) \le \\
C_1\left\{\gamma^{-1}\int_{|x|\ge\frac{R}{2}}|g(x)|^2dx + \frac{\nu E}{R^2} +
\frac{U}{R}\|\omega(\cdot, s)\|_{L^4({\mathbb R}^2)}^2\right\}
\ea
$$ 
where $U$ is a time independent bound on $\|u\|_{L^2({\mathbb R}^2)}$, depending only on $\gamma$, $\|u_0\|_{L^2({\mathbb R}^2)}$ and on $\|f\|_{L^2({\mathbb R}^2)}$, and 
$E$ is a time independent bound on the enstrophy, depending only
on $\gamma$, $\|\omega_0\|_{L^2({\mathbb R}^2)}$ and $\|g\|_{L^2({\mathbb R}^2)}$. We observe that
$$
\nu \int_0^{\infty}e^{\gamma (t-s)}\|\nabla\omega(\cdot ,s)\|_{L^2({\mathbb R}^2)}^2ds
$$
is bounded in terms of $\gamma$, initial enstrophy and the norm of $g$ in $L^2({\mathbb R}^2)$, a fact that follows immediately from the enstrophy balance. From the uniform bound on enstrophy and a Sobolev embedding theorem we deduce that
$$
\int_0^{\infty}e^{\gamma (t-s)}\|\omega(\cdot, s)\|_{L^4({\mathbb R}^2)}^2ds \le F
$$
where $F$ is bounded in terms of $\gamma$,  the viscosity, initial enstrophy and norm of $g$ in $L^2({\mathbb R}^2)$. It then follows that 
$$
Y_R(t) \le e^{-\gamma t}Y_R(0) + C_2\left\{\gamma^{-1}\int_{|x|\ge\frac{R}{2}}|g(x)|^2dx + \frac{E\nu}{\gamma R^2} + \frac{UF}{R}\right\}
$$
Choosing $R$ large enough proves the claim. We note 
that $R$ can be chosen uniformly for 
all initial vorticities $\omega_0\in L^2({\mathbb R}^2)$ that are uniformly bounded in $L^2({\mathbb R}^2)$ and  satisfy a uniform centering property (see below).

We are going to use the notation $\langle f, g\rangle = \int_{{\mathbb R}^2}f(x)g(x)dx$, and sometimes write $S^{NS,\gamma}(t)(\omega_0)$
for the vorticity $\omega(x,t)$ solution of (\ref{omeq}).
\begin{thm}\label{compindi}
Let $\omega_0\in X$ where $X \subset L^{2}({\mathbb R}^2)$ is a bounded set that satisfies the uniform centering property $\forall\epsilon>0$, $\exists R>0$, $\forall\omega_0\in X$
$$
\int_{|x|\ge R}\left |\omega_0(x)\right|^2dx \le \epsilon.
$$ 
Then, for any $t_0>0$, the set
$$
O^{+}(t_0, X) = cl\left \{ S^{NS,\gamma}(t)\omega_0\,\, \left |\, \omega_0\in X, \, t\ge t_0\right. \right\}
$$
(where $cl(O)$ is the $L^2({\mathbb R}^2)$ closure of the set $O$) is compact in $L^2({\mathbb R}^2)$.
\end{thm}
The proof of this theorem follows from an uniform bound in $H^1({\mathbb R}^2)$ for $\omega (t)$ for $t\ge t_0$ and the uniform ``no-travel''property of the previous theorem.

\section{Stationary Deterministic Solutions.}
Let $(u^{(\nu)}, \omega^{(\nu)})$ be a sequence of solutions of
\be
\left \{
\ba
-\nu\Delta u + \gamma u + \nabla p + u\cdot\nabla u  = f,\\
\nabla\cdot u = 0
\ea
\right.
\la{ust}
\ee
and the corresponding vorticity equation
\be
\left \{
\ba
\gamma \omega + u\cdot\nabla\omega -\nu\Delta \omega = g,\\
\omega = \nabla^{\perp} \cdot u.
\ea
\right .
\la{stomeq}
\ee 
We let $\nu \to 0$ but keep $f, g, \gamma$ fixed. The solutions $u^{(\nu)}$ exist, are smooth and decay rapidly at infinity. Moreover, the energy balance
$$
\gamma \|u^{(\nu)}\|_{L^2({{\mathbb R}^2})}^2 + \nu\|\nabla u^{(\nu)}\|^2_{L^2({\mathbb R}^2)} =
\int_{{\mathbb R}^2} f\cdot u^{(\nu)}dx
$$
implies that the sequence $u^{(\nu)}$ is bounded in $L^2(\mathbb{R}^2)$. 

The enstrophy balance
\be
\gamma \|\omega\|^2_{L^2({\mathbb R}^2)} + \nu \|\nabla \omega\|_{L^2({\mathbb{R}^2})}^2 = \int_{{\mathbb R}^2} g\omega dx
\la{viscensbal}
\ee
implies that the sequence $\omega^{(\nu)}$ is bounded in $L^2({\mathbb R}^2)$. 
 Passing to a subsequence, we consider the weak limit 
\be
\omega^{(0)} = w-\lim_{\nu\to 0}\omega^{(\nu)}
\la{omegz}
\ee
in $L^2(\mathbb{R}^2)$.
Because of the compact embedding $H^1(\mathbb{R}^2)^2\subset L^q(\Omega)^2$, for any relatively compact open set $\Omega \subset {\mathbb R}^2$ and any $1\le q<\infty$,  we may assume, by passing to a subsequence, that $u^{(\nu)}=K\star\omega^{(\nu)}$ (where $K = \frac{1}{2\pi}\frac{x^{\perp}}{|x|^2}$ is the Biot-Savart kernel) converge to $u^{(0)}$ strongly in $L^q(\Omega)^2$.
\begin{thm} The function $\omega^{(0)}$ is a renormalized solution of the inviscid equation
\be
\left\{
\ba
\gamma \omega^{(0)} + u^{(0)}\cdot\nabla \omega^{(0)} = g\\
\omega^{(0)} = \nabla^{\perp}\cdot u^{(0)}
\ea
\right.
\la{reno}
\ee
In addition, $\omega^{(0)}\in L^2({\mathbb R}^2)$, $u^{(0)}\in H^1({\mathbb R}^2)$,
the equation holds in $W^{-1, q}_{loc}({\mathbb R}^2)$ for any  $1<q<2$, and
\be
\gamma \|\omega^{(0)}\|^2_{L^2({\mathbb R}^2)} = \int_{{\mathbb{R}^2}} g\omega^{(0)}dx
\la{bala}
\ee

holds.
\end{thm}
\noindent{\bf Remark.} Renormalized solutions have been introduced in 
(\cite{dpl}). The existence of weak solutions for damped and driven Euler equations
using a vanishing viscosity method was obtained in (\cite{bct}).

\noindent{\bf Proof.} 

The facts that $\omega^{(0)}\in L^2({\mathbb R}^2)$, and $u^{(0)}\in L^2({\mathbb R}^2)^2$ follow from the construction and uniform bounds on the solutions $u^{(\nu)}, \omega^{(\nu)}$. 
Furthermore, the solutions $u^{(\nu)}$ are bounded in $H^1(\mathbb{R}^2)^{2}$ and converge strongly in $L^2_{loc}(\mathbb{R}^2)^{2}$ to $u^{(0)}$. The vorticities are bounded in $L^2(\mathbb{R}^2)$ and converge weakly. If $\phi$ is a test function then  $(u^{(\nu)}\cdot\nabla)\phi$ converge strongly to $(u^{(0)}\cdot\nabla)\phi$ in $L^2(\mathbb{R}^2)^2$ and, because the scalar product of weak and strong convergent sequences is convergent, we have:
$$
\lim_{\nu\to 0}\int_{{\mathbb R}^2}u^{(\nu)}\cdot\nabla\phi \omega^{(\nu)} dx =
\int_{{\mathbb R}^2}u^{(0)}\cdot\nabla\phi\omega^{(0)}dx.
$$
This means that $u^{(0)}, \omega^{(0)}$ is a weak solution of the inviscid
equation. Because $u^{(0)}\in  H^{1}(\mathbb{R}^2)^2$, $\omega^{(0)}\in L^{2}(\mathbb{R}^2)$ and $g\in L^{2}(\mathbb{R}^2)$, we are under the conditions of consistency in (\cite{dpl}), Thm II. 3, and the same proof applied 
to our case shows  that $u^0$, $w^0$ is a renormalized solution of the inviscid equation, that is,

\begin{equation}
\label{reuler}
\gamma \omega^{(0)}\beta'(\omega^{(0)}) + u^{(0)}\cdot\nabla \beta(\omega^{(0)}) = g\beta'(\omega^{(0)})
\end{equation} 
holds in the sense of distributions for any $\beta\in C^1$ that is bounded, has bounded derivative and vanishes near the origin. We present the proof here,
for the sake of completeness.
It is easy to prove (see Lemma II$.1$ in \cite{dpl}), that if $u^0\in (W^{1,2}_{loc}(\mathbb{R}^2))^2$, $\omega^0\in L_{loc}^2(\mathbb{R}^2)$, then
\be
\label{dplconv}
\left(u^0\cdot\nabla\omega^0\right)\star j_{\epsilon}-u^0\cdot\nabla\left(\omega^0\star j_{\epsilon}\right)\rightarrow 0\quad \text{in}\; L^1_{loc}(\mathbb{R}^2)
\ee
when $\epsilon$ tends to zero. Here (and hereafter)  
$j_{\epsilon}$ is a standard mollifier -- $j_{\epsilon}(z) = \epsilon^{-2}j(\epsilon^{-1}z)$ with $j(z)$ a fixed smooth, even, compactly supported nonnegative function with $\int j(z)dz = 1$ --  and $a\star b$ denotes convolution.

Then, considering the mollified functions
$\omega^0_\epsilon=\omega \star j_\epsilon$, $u^0_\epsilon =u^0\star j_\epsilon$ and $g_\epsilon=g\star j_\epsilon$, it follows immediately from (\ref{dplconv}) that
\be
\label{dplconveq}
u^0\cdot\nabla\omega^0_{\epsilon} + \gamma\omega^0_\epsilon-g_\epsilon=q_\epsilon,
\ee 
holds in the sense of distributions, and $q_\epsilon$ converges to zero in $L^1_{loc}(\mathbb{R}^2)$ as $\epsilon$ tends to zero. From this equation, we obtain that if $\beta\in C^1(\mathbb{R})$, and $\beta$ is bounded with bounded first derivative, then  

\be
\label{dplrenormeq}
u^0\cdot\nabla\beta(\omega^0_{\epsilon}) + \gamma\omega^0_\epsilon\beta'(\omega^0_{\epsilon})-g_\epsilon\beta'(\omega^0_{\epsilon}) = q_\epsilon\beta'(\omega^0_{\epsilon}).
\ee 
also holds in the sense of distributions. Letting $\epsilon$ tend to zero, we prove (\ref{reuler}).

In order to prove $(\ref{bala})$, we mollify $b = \beta(\omega^{(0)})$, where $\beta$ is a $C^1$ function with compact support
$$
b_{\epsilon} = b \star j_{\epsilon}.
$$
We use the identity (\cite{cet}) 
\begin{equation}
(u \otimes b)_{\epsilon}-u_{\epsilon}\otimes b_{\epsilon} = 
\rho_{\epsilon}(u,b),
\label{flux}
\end{equation}
with
\begin{equation}
\rho_{\epsilon}(u,b) =
r_{\epsilon}(u, b) -
(u-u_{\epsilon})\otimes(b- b_{\epsilon}),
\label{rho}
\end{equation} 
and with
$$
r_{\epsilon} (u, b) = \int_{{\mathbb R}^2}j(z)(u(x-\epsilon z)-u(x))\otimes (b(x-\epsilon z)-b(x))dz.
$$
Because 
$$
\int_{\mathbb{R}^2}\text{Tr}\left(\left(u_{\epsilon}\otimes b_{\epsilon}\right)\nabla b_{\epsilon}\right)dx=0,
$$
it follows that
\be
\int_{{\mathbb R}^2}(b u)_{\epsilon}\cdot \nabla b_{\epsilon}dx =
\int_{{\mathbb R}^2}r_{\epsilon}(u, b)\cdot\nabla b_{\epsilon} -
\int_{{\mathbb R}^2}(u-u_{\epsilon})(b -b_{\epsilon})\nabla b_{\epsilon}dx.
\la{molbal}
\ee
Now  $b= \beta(\omega)\in L^1\cap L^{\infty}({\mathbb R}^2)$ and we can pass to the limit in  $(\ref{molbal})$ using the fact that $u-u_{\epsilon}$ is $O(\epsilon)$ in $L^2({\mathbb R}^2)$ (because of the uniform bound in $H^1({\mathbb R}^2)$), and working in $L^4({\mathbb R}^2)$ with $b$ : $\nabla b_{\epsilon}$ is $O(\epsilon)^{-1}$ in $L^4({\mathbb R}^2)$, and $b-b_{\epsilon}$ converges to zero in $L^4({\mathbb R}^2)$. We deduce that
$$
\gamma\int_{{\mathbb R}^2}\omega^{(0)}\beta'(\omega^{(0)})\beta(\omega^{(0)})dx  = \int_{{\mathbb R}^2}g\beta'(\omega^{(0)})\beta(\omega^{(0)})dx
$$
holds for any $\beta\in C^1$ with compact support.
Taking a sequence of functions that approximate $\beta(\omega) = \omega$, with $\beta'$ uniformly bounded, we deduce (\ref{bala}).

\begin{thm} Let $u^{(\nu)}, \omega^{(\nu)}$ be a sequence of solutions of (\ref{ust}, \ref{stomeq}).
Then the enstrophy dissipation vanishes in the limit $\nu\to 0$:
$$
\lim_{\nu\to 0}\nu\int_{{\mathbb R}^2}|\nabla\omega^{(\nu)}|^2dx = 0
$$
holds.
\end{thm}

\noindent{\bf{Proof.}} 
Taking  the limit superior in the enstrophy balance equation (\ref{viscensbal}),  using Fatou's lemma and the fact that $\omega^{(\nu)}$ converge to $\omega^{(0)}$ weakly in $L^{2}(\mathbb{R}^2)$, we have: 
\be
\begin{aligned}
\limsup_{\nu\rightarrow 0}\nu\|\nabla \omega^{(\nu)}\|_{L^2({\mathbb{R}^2})}^2&\leq\limsup_{\nu\rightarrow 0}\int_{{\mathbb R}^2} g\omega^{(\nu)} dx-  \liminf_{\nu\rightarrow 0} \gamma\|\omega^{(\nu)}\|^2_{L^2({\mathbb R}^2)}  \\
&\leq\int_{{\mathbb R}^2} g\omega^{(0)} dx-\gamma \|\omega^{(0)}\|^2_{L^2({\mathbb R}^2)} =0.
\end{aligned}
\ee

\section{Stationary Statistical Solutions}

In this section we follow the methods of Foias, see \cite{foias72, foias73, fmrt}, and define a notion of  a stationary statistical solution of the damped and driven incompressible Navier-Stokes equations in the vorticity phase space. 
The solution  is a Borel probability measure in $L^2(\mathbb{R}^2)$.
 We note that $L^2(\mathbb{R}^2)$ is a separable Hilbert space and therefore the Borel $\sigma$-algebra associated to the strong (norm) topology is the same as the Borel $\sigma$-algebra associated to the weak topology.  (Indeed, any open set is a countable union of open balls, any open ball is a countable union of closed balls and closed balls are convex, hence weakly closed, hence weakly Borel.)

\begin{sssvort}
\label{sssvort}

A stationary statistical solution of the damped  and driven \\
Navier-Stokes equation (SSSNS) in vorticity phase space is a Borel probability measure $\mu^{\nu}$ in $L^2(\mathbb{R}^2)$ such that
\begin{enumerate}

\item[(1)]
\label{sssv1}
$\displaystyle{\int_{L^2(\mathbb{R}^2)}\left\|\omega\right\|_{H^1(\mathbb{R}^2)}^{2}d\mu^{\nu}(\omega)<\infty ,}$

\item[(2)]
\label{sssv2}
$\displaystyle{\int_{L^2(\mathbb{R}^2)}\langle u \cdot\nabla\omega +\gamma\omega-g,\Psi'(\omega)\rangle+\nu\langle\nabla_x\omega,\nabla_x\Psi'(\omega)\rangle  d\mu^{\nu}(\omega)=0}
\\  \text{for any test functional}\; \Psi\in\mathcal{T},\;\, 
\text{with }\; \, u = \frac{1}{2\pi}\frac{x^{\perp}}{|x|^2}\star \omega,\,\,
\text{and} $     

\item[(3)]
\label{sssv3}
$
\displaystyle{\int_{ E_1\leq\left\|\omega\right\|_{L^2({\mathbb R}^2)}\leq E_2}\left\{\gamma\left\|\omega\right\|_{L^2(\mathbb{R}^2)}^2+\nu\left\|\omega\right\|_{H^1(\mathbb{R}^2)}^{2}-\langle\mathbf{g},\omega\rangle\right\}d\mu^{\nu}(\omega)\le 0},\\
\quad E_1, E_2 \ge 0$.

\end{enumerate}
\end{sssvort}
The class of cylindrical test functions ${\mathcal{T}}$ is given by:

\begin{testII}
\label{testII}
The class of test functions $\mathcal{T}$ is the set of functions $\Psi:L^2(\mathbb{R}^2)\rightarrow\mathbb{R}$ of the form
\be
\Psi(\omega):=\Psi_I(\omega)=\psi\left(\langle\omega,\mathbf{w}_{1}\rangle,\ldots,\langle\omega,\mathbf{w}_{m}\rangle\right),
\ee
or
\be
\Psi(\omega):=\Psi_{\epsilon}(\omega)=\psi\left(\langle\alpha_{\epsilon}(\omega),\mathbf{w}_{1}\rangle,\ldots,\langle\alpha_{\epsilon}(\omega),\mathbf{w}_{m}\rangle\right),
\ee
where $\psi$ is a $C^{1}$ scalar valued function defined on $\mathbb{R}^{m}$,	
$m\in\mathbb{N}$;
$\mathbf{w}_{1},\ldots,\mathbf{w}_{m}$ belong to $C^2_0({\mathbb R}^2)$ and
$$
\alpha_{\epsilon}(\omega)= J_{\epsilon}\beta(J_{\epsilon}\omega),
$$
where $\beta\in C^3$ is a compactly supported function of one real variable,
and $J_\epsilon$ is the convolution operator
$$
J_\epsilon(\omega ) = j_{\epsilon}\star\omega.
$$
with $j\ge 0$ a fixed smooth, nonnegative, even  $(j(-z) = j(z))$ function
supported in $|z|\le 1$ and with $\int_{\mathbb{R}^2}j(z)dz =1$.
\end{testII}

The test functions $\Psi$ used  in the definition are all 
locally bounded and weakly sequentially continuous in $L^2({\mathbb R}^2)$. We note the trivial but very important distinction between weakly continuous and strongly continuous functions defined on $L^2({\mathbb R}^2)$: any weakly continuous function is strongly continuous, but there exist strongly continuous functions -- for instance, the norm -  that are not weakly continuous.  Because the SSSNS is a Borel probability, bounded continuous functions are integrable. In the sequel we will pass to weak limits of SSSNS, $\mu^{\nu}\to\mu^{E}$, and then the distinction between weakly continuous
and continuous functions is important: although for strongly continuous functions $\Psi$ the integrals $\int\Psi d\mu^{\nu}$ are defined and finite, it is only for weakly continuous functions $\Psi$ that the weak limit $\lim_{\nu\to 0}\int\Psi d\mu^{\nu} = \int\Psi d\mu^{E}$ holds by definition. We will obtain stronger information as well, but that needs to be proved carefully.

We discuss now the definition of SSSNS and comment on
its mathematical soundness. We will also verify the fact that for each test function, the  integrand in (2) is a weakly continuous function on $L^2$. We start by making sense of (1): the integrand can be viewed as a Borel measurable function defined for all $\omega\in L^2({\mathbb R}^2)$, equal to infinity for $\omega\not\in H^{1}({\mathbb R}^2)$. The fact that this function is Borel measurable follows from the fact that $\|\omega\|^2_{H^1}$ is everywhere the limit of the sequence of continuous (hence measurable) functions $\|J_{\epsilon}\omega\|^2_{H^1}$ 
obtained by taking a fixed a sequence $\epsilon\to 0$ and convolving with a mollifier. The requirement (3) is a local enstrophy balance; it implies (\cite{fmrt}) that the SSSNS has bounded support. We define the set   
\be
\label{Bdef}
B:=\left\{\omega\in L^2(\mathbb{R}^2); \left\|\omega\right\|_{L^2({\mathbb R}^2)}\leq \frac{\left\|g\right\|_{L^2({\mathbb R}^2)}}{\gamma}\right\}.
\ee

\begin{boundsupp}
The support of any stationary statistical solution of the damped and
driven Navier-Stokes equations in vorticity phase space
is included in the bounded set in $B\subset L^2(\mathbb{R}^2)$:
\begin{equation}
\label{boundsupp}
\supp\mu^{\nu}\subset B.
\end{equation}
\end{boundsupp}

\noindent{\bf{Proof.}}
It follows from Definition \ref{sssvort} item (3), that if 
\[
E=\left\{\omega\in L^2(\mathbb{R}^2); E_1^2\leq\left\|\omega\right\|_{L^2({\mathbb R}^2)}^2\leq E_2^2\right\},
\] 
then
$$
\ba
\gamma\int_E\left\|\omega\right\|_{L^2({\mathbb R}^2)}^2d\mu^{\nu}(\omega)\leq\left\|g\right\|_{L^2({\mathbb R}^2)}\int_E\left\|\omega\right\|_{L^2({\mathbb R}^2)}d\mu^{\nu}(\omega)
\\
\leq\left\|g\right\|_{L^2({\mathbb R}^2)}\left(\int_E\left\|\omega\right\|_{L^2({\mathbb R}^2)}^2d\mu^{\nu}(\omega)\right)^{1/2}.
\ea
$$
Hence,
\[
\int_E\left\|\omega\right\|_{L^2({\mathbb R}^2)}^2d\mu^{\nu}(\omega)\leq\frac{\left\|g\right\|_{L^2({\mathbb R}^2)}^2}{\gamma^2}.
\]
Thus,
\begin{equation}
\label{suppargument}
\int_E\left(\left\|\omega\right\|_{L^2({\mathbb R}^2)}^2-\frac{\left\|g\right\|_{L^2({\mathbb R}^2)}^2}{\gamma^2}\right)d\mu^{\nu}(\omega)\leq 0.
\end{equation}
If ${E}_1^2=\left\|g\right\|_{L^2({\mathbb R}^2)}^2/\gamma^2$ and ${E}_2\to\infty$, then by \eqref{suppargument}, we have $\mu({E})=0$, and the result follows immediately.$\quad\Box$

We compute now $\Psi'$ for the test functions $\Psi\in {\mathcal T}$.
Clearly  $\alpha_{\epsilon}: L^{2}({\mathbb R}^2)\mapsto L^{2}({\mathbb R}^2)$  continuously differentiable and bounded uniformly on bounded sets of $L^2({\mathbb R}^2)$; moreover
\be
\alpha_{\epsilon}'(\omega)\cdot \phi = ((\beta'(\omega_\epsilon))_\epsilon \phi_\epsilon)_\epsilon,\;\forall \phi\in L^1_{loc}(\mathbb{R}^2).
\la{alphaprime}
\ee
For $\Psi_{\epsilon}\in{\mathcal T}$ we have thus
\be
\label{Iprod}
\nabla_\omega\Psi_{\epsilon}(\omega)\cdot \phi =\sum_{j=1}^{m}\partial_{j}\psi(\langle\alpha_{\epsilon}(\omega),\mathbf{w}_{1}\rangle,\ldots,\langle\alpha_{\epsilon(\omega)},\mathbf{w}_{m}\rangle)\langle(\beta'(\omega_\epsilon) \phi_\epsilon)_\epsilon,\mathbf{w}_{j}\rangle,
\ee
and for $\psi_I\in{\mathcal T}$ we have
\be
\label{epprod}
\nabla_\omega\Psi_I(\omega)\cdot \phi=\sum_{j=1}^{m}\partial_{j}\psi(\langle\omega,\mathbf{w}_{1}\rangle,\ldots,\langle\omega,\mathbf{w}_{m})\rangle)\langle \phi,\mathbf{w}_{j}\rangle,
\ee
where $\partial_{j}\psi$ denotes the derivative of $\psi$ with
respect to its $j$-th variable.
Clearly, in both cases, $\phi \mapsto \nabla_\omega\Psi(\omega)\cdot \phi$ is a bounded linear continuous functional on $L^2({\mathbb R}^2)$ and thus, by the Riesz representation theorem, there exists an element $\Psi'(\omega)\in L^2({\mathbb R}^2)$ such that
\[
\nabla_\omega\Psi(\omega)\cdot v=\langle\Psi'(\omega),v\rangle,\quad\forall v\in L^2({\mathbb R}^2).
\]
This is the identification implied in the shorthand notation $\Psi'(\omega)$ used in Definition $\ref{sssvort}$. For instance
\be
\Psi_I'(\omega)=\sum_{j=1}^{m}\partial_{j}\psi(\langle\omega,\mathbf{w}_{1}\rangle,\ldots,\langle\omega,\mathbf{w}_{m}\rangle)\mathbf{w}_{j}.
\la{psiprimei}
\ee
Consequently
\be
\label{Ibinfty}
\partial_x^{(k)}\Psi_I'(\omega)=\sum_{j=1}^{m}\partial_{j}\psi(\langle\omega,\mathbf{w}_{1}\rangle,\ldots,\langle\omega,\mathbf{w}_{m}\rangle)\partial_x^{(k)}\mathbf{w}_{j}.
\ee
for any multi-index $k$ with $|k|\le 2$. For $\Psi_{\epsilon}$, a similar computation yields 
\be
\partial_k^{(k)}\Psi_{\epsilon}'(\omega) = \sum_{j=1}^m\partial_j\psi(\langle
\alpha_{\epsilon}(\omega), {\mathbf{w}}_1\rangle, \ldots, \langle\alpha_{\epsilon}(\omega), {\mathbf{w}}_m\rangle)\partial_x^{(k)}
\left (\beta'(\omega_{\epsilon}){{\mathbf{w}}_j}_{\epsilon}\right)_{\epsilon}.
\la{psiprimee}
\ee

\begin{lemma} \label{psiregs} Let $\Psi\in {\mathcal{T}}$ and $\omega\in B\subset L^2({\mathbb {R}}^2)$ with $B$ a bounded set in $L^{2}({\mathbb {R}}^2)$. Then $\Psi'(\omega)\in C_0^2({\mathbb R}^2)$, and there exists a constant depending only on  $\Psi$ and $B$ such that
\be
\|\Psi'(\omega)\|_{W^{2,2}({\mathbb R}^2)} + \|\Psi'(\omega)\|_{W^{2,\infty}({\mathbb R}^2)} \le C
\la{regu}
\ee
holds for all $\omega \in B$.

Consider, for any 
$$
F_i : L^{2}({\mathbb R}^2) \to {\mathbb R}
$$
$i=1, 2, 3$ given by
$$
F_1(\omega) = \langle \Psi'(\omega),\gamma\omega -g\rangle,
$$
$$
F_2(\omega) = \langle \nax\Psi'(\omega), \nax \omega\rangle
$$
and
$$
F_3(\omega) = \langle \Psi'(\omega), u\cdot\nabla\omega \rangle,\,\quad  u = \frac{1}{2\pi}\frac{x^{\perp}}{|x|^2}\star \omega.
$$
These three maps are well defined for $\omega\in L^{2}({\mathbb R}^2)$, weakly continuous and bounded uniformly on bounded sets $B\subset L^{2}({\mathbb R}^2)$.
\end{lemma}
\noindent{\bf{Remarks.}} If $\beta\in C_{0}^{k+1}({\mathbb{R}})$ and
${\mathbf{w}}_j\in C_0^{k}({\mathbb R}^2)$ then $\Psi'(\omega)\in C_0^{k}({\mathbb{R}}^2)$. The expressions $\langle\nax\Psi'(\omega), \nax\omega\rangle = 
-\langle \Delta_x\Psi'(\omega), \omega\rangle$ and $\langle \Psi'(\omega), u\cdot\nabla\omega\rangle = -\langle u\cdot\nax\Psi'(\omega), \omega\rangle$ make sense for
$k\ge 2$, $u\in L^2_{loc}({\mathbb R}^2)$, $\omega\in B$.

\noindent{\bf{Proof.}} 
It is easy to see that
$ \partial_{x}^{(k)}(\Psi_{\epsilon})'(\omega)$ and $\partial_{x}^{(k)}\Psi_I'(\omega)$ are uniformly bounded in $L^\infty({\mathbb R}^2)\cap L^2({\mathbb R}^2)$ for all $\omega\in B$, $|k|\le 2$. This is verified for $\partial_{x}^{(k)}\Psi_I'(\omega)$ directly by inspection of  \eqref{Ibinfty} and for
 $\partial_{x}^{(k)}\Psi_{\epsilon}'$ by inspection of (\ref{psiprimee}).
We check the bounds for $\Psi_{\epsilon}$:
As $\alpha^\epsilon(\omega)$ is bounded on bounded subsets of $L^2(\mathbb{R}^2)$, and $\psi$ is of class $C^1$, we have that
\be
\label{psibound}
\left|\partial_{j}\psi(\langle\alpha^\epsilon(\omega),\mathbf{w}_{1}\rangle,\ldots,\langle\alpha^\epsilon(\omega),\mathbf{w}_{m}\rangle)\right|\leq C,\quad\forall\omega\in B.
\ee 
The fact that $\beta'(\omega_{\epsilon}) \in L^{\infty}({\mathbb R}^2)$ is bounded uniformly for $\omega\in B$ implies that 
\be
\label{difflinn}
\left\|\partial_x^{(k)}\left ((\beta'(\omega_{\epsilon}) {{\mathbf{w}}_j}_{\epsilon}\right )_{\epsilon}\right\|_{L^p({\mathbb{R}}^2)}\leq \frac{C}{\epsilon^{|k|}}\left\| {\mathbf{w}}_j\right\|_{L^p({\mathbb{R}}^2)}
\ee
holds uniformly, for all $p$, $1\le p\le\infty$. By \eqref{psibound} and \eqref{difflinn}, we have from (\ref{psiprimee}) that
\be
\label{difflin2}
\left\|\partial_x^{(k)}\Psi_{\epsilon}'(\omega)\right\|_{L^p({\mathbb{R}}^2)}
\le \frac{C_p}{\epsilon^{|k|}}
\ee
holds for $1\le p\le\infty$ with $C_p$ uniform all $\omega\in B$.
Thus, $\partial_{x}^{(k)}\Psi'(\omega)$ are bounded in  $L^\infty({\mathbb R}^2)\cap L^2({\mathbb R}^2)$.

Concerning the statements about the maps $F_i$, we start with
\be
\label{F1def}
F_1(\omega)=\langle\Psi'(\omega),\gamma\omega-g\rangle=\nabla_\omega\Psi(\omega)\cdot (\gamma\omega-g);\quad \omega\in L^2(\mathbb{R}^2).
\ee
This function is weakly continuous. Indeed, for $\Psi_I$ we have by (\ref{psiprimei}) 
\[
\begin{aligned}
& \langle\Psi_I'(\omega),\gamma\omega-g\rangle \\ & = \sum_{j=1}^{m} \partial_{j}\psi(\langle\omega,\mathbf{w}_{1}\rangle,\ldots,\langle\omega,\mathbf{w}_{m}\rangle)\langle\mathbf{w}_{j},\gamma\omega-g\rangle.
\end{aligned}
\]
and it is clear that this is a weakly continuous function of
$\omega\in L^2({\mathbb R}^2)$. It is also quite obvious that it is uniformly bounded for $\omega \in B$. In the case of $\Psi_{\epsilon}$, by \eqref{Iprod} we have
\[
\begin{aligned}
&\nabla_\omega\Psi_\epsilon(\omega)\cdot (\gamma\omega-g)\\
& = \sum_{j=1}^{m} \partial_{j}\psi(\langle\alpha^\epsilon(\omega),\mathbf{w}_{1}\rangle,\ldots,\langle\alpha^\epsilon(\omega),\mathbf{w}_{m}\rangle)\langle((\beta'(\omega_{\epsilon})(\gamma\omega-g)_{\epsilon})_{\epsilon},\mathbf{w}_{j}\rangle.
\end{aligned}
\]
The  weak continuity here follows from the fact that if $\omega^j$ converges weakly to $\omega$ then $\omega^j_{\epsilon} \to \omega_{\epsilon}$ converge pointwise, and it is bounded. Consequently, $(\beta'(\omega^j_{\epsilon})(\gamma\omega^j-g)_{\epsilon})_{\epsilon}$ converges pointwise and is uniformly bounded. Therefore  we can use the Lebesgue dominated convergence theorem in the integral against a fixed ${\mathbf{w}}$ from the finite list ${\mathbf{w}}_1, \dots {\mathbf{w}}_m$ appearing in $\Psi^{\epsilon}$. It is also clear that
\begin{equation}
\left\|\alpha_{\epsilon}'(\omega)\cdot (\gamma\omega-g)\right\|_{L^2({\mathbb R}^2)}\leq C(\left\|\omega\right\|_{L^2({\mathbb R}^2)}+\left\|g\right\|_{L^2({\mathbb R}^2)}),\quad\forall  \omega\in L^2(\mathbb{R}^2). 
\end{equation}
Thus, we have
\be
\label{statomegabound}
F_1(\omega)\leq c_3(\left\|\omega\right\|_{L^2({\mathbb R}^2)}+\left\|g\right\|_{L^2({\mathbb R}^2)})\leq C;\quad\forall  \omega\in B.
\ee
Therefore, $F_1(\omega)$ is weakly continuous and bounded uniformly for $\omega\in B$. 

The fact that $F_2$ is well defined follows from the fact that $\Delta \Psi'(\omega)\in L^2({\mathbb R}^2)$ and
$$
F_2(\omega) = -\langle \Delta\Psi'(\omega), \omega\rangle.
$$
The weak continuity for
$F_2$ follows as for $F_1$: in the case of $\Psi_I$ it is straightforward, and in the case of $\Psi_{\epsilon}$ it follows because weak convergence becomes pointwise convergence and we can apply the Lebesgue dominated convergence theorem. 

For $F_3$, we note first that, if  $u = \frac{1}{2\pi}\frac{x^{\perp}}{|x|^2}\star \omega$, then, by classical singular integral theory (\cite{stein}) $u\in L^r_{loc}({\mathbb R}^2)$,
$r<\infty$, and $\nabla\cdot u =0$. Because $\nax\Psi'(\omega)$ is bounded and compactly supported, $u\cdot\nax\Psi'(\omega)\in L^2({\mathbb R}^2)$ and
$$
F_3(\omega) = -\langle u\cdot\nax\Psi'(\omega), \omega\rangle
$$
is well defined.

If $\omega_k$ converge weakly in $L^{2}({\mathbb R}^2)$ to $\omega$, then the corresponding velocities $u_k$ converge strongly to $u$ in $L^2$ on compact sets$K$, by the compact embedding $H^1(K)\subset\subset L^2(K)$. The case of $\Psi_I$ follows then because the functions ${\mathbf{w}}_j$ in the list ${\mathbf{w}}_1, \dots,  {\mathbf{w}}_m$ have compact supports, and therefore the functions
$u_k\cdot\nabla {\mathbf{w}}_j$ converge strongly as $k\to\infty$ to 
$u\cdot\nabla {\mathbf{w}}_j$ in $L^2({\mathbb R}^2)$. The scalar products 
$$
\langle u_k\cdot\nabla \omega_k , {\mathbf{w}}_j\rangle = -\langle \omega_k, u_k\cdot\nabla {\mathbf{w}}_j\rangle
$$
converge, as $k\to\infty$  to $-\langle \omega, u\cdot\nabla {\mathbf{w}}_j\rangle$,
because the scalar products of weakly convergent and strongly convergent sequences converge. Therefore the function $F_3$ is weakly continuous for this class of test functions. It is easy to see that the function is uniformly bounded locally in $L^2({\mathbb R}^2)$. In the case of $\Psi_{\epsilon}$, a similar argument shows that, if $\omega_k$ converges weakly in $L^2({\mathbb R}^2)$ to $\omega$, then
$$
\left (u_k\cdot\nabla \omega_k\right)_\epsilon (x)\to \left (u\cdot\nabla \omega\right )_{\epsilon}(x)
$$
holds for each $x\in {\mathbb R}^2$, and these functions are uniformly bounded
as $x\in{\mathbb R}^2$. Also, the functions $\beta'((\omega_k)_{\epsilon})$
converge pointwise and are bounded. This implies that $F_3$ is weakly continuous; the uniform boundedness is easily verified. $\quad\Box$

We define the notion of renormalized stationary statistical solution of the Euler equation.
\begin{Eulerdef}
\label{Eulerdef}
A Borel probability measure $\mu^0$ on $\mathbb{R}^2$ is a renormalized stationary statistical solution of the damped and driven Euler equation if 
\be
\int_{L^2(\mathbb{R}^2)}\langle u\cdot\nabla \omega + \gamma\omega-g,\Psi'(\omega)\rangle d\mu^{0}(\omega) = 0
\la{renor}
\ee
(with $u = \frac{1}{2\pi}\frac{x^{\perp}}{|x|^2}\star \omega$) holds for any test functional  $\Psi\in\mathcal{T}$.

We say that a renormalized stationary statistical solution $\mu^{0}$ 
of the Euler equation satisfies the enstrophy balance if
\be
\int_{L^2(\mathbb{R}^2)}\left\{\gamma\left\|\omega\right\|_{L^2(\mathbb{R}^2)}^2-\langle\mathbf{g},\omega\rangle\right\}d\mu^{0}(\omega)=
  0
\la{bal}
\ee
holds.
\end{Eulerdef}

We recall  Prokhorov's theorem (see for instance \cite{smolfomin}):
\begin{Prokhorov}
\label{Prokhorov}
Let $X$ be a complete separable metrizable topological space, and let $\mathcal{M}$ be a set of Borel probability measures on $X$. For each sequence in $\mathcal{M}$ to contain a weakly convergent subsequence it is sufficient that for each $\epsilon >0$, there is a compact subset $K$ of $X$ such that $\mu(X\setminus K)<\epsilon$ for each $\mu\in\mathcal{M}$. 
\end{Prokhorov}
We recall that a sequence of Borel probability measures $\pi_n$ on a topological space $X$ converges weakly to a Borel probability measure $\pi$ on $X$ if for every continuous bounded real-valued function $\Psi$ on $X$
\be
\lim_{n\to\infty}\int_{X}\Psi(s)d\pi_n(s) = \int_{X}\Psi(s)d\pi(s).
\ee

\begin{ilimit}
\label{ilimit} 
Given a sequence of stationary statistical solutions of the damped and 
driven NSE in vorticity phase space, $\left\{\mu^\nu\right\}$, 
with $\nu\rightarrow 0$, there exists a subsequence, denoted also
$\left\{\mu^\nu\right\}$, and a Borel probability measure $\mu^0$ on  $L^2(\mathbb{R}^2)$, such that
\begin{equation}
\label{intlimit}
\lim_{\nu\to 0}\int_{L^2(\mathbb{R}^2)}\Phi(\omega)d\mu^{\nu}(\omega)= \int_{L^2(\mathbb{R}^2)}\Phi(\omega)d\mu^{0}(\omega),
\end{equation}
holds for all weakly continuous, locally  bounded real-valued functions $\Phi$. Furthermore, the weak limit measure $\mu^0$ is a renormalized stationary statistical solution of the damped and driven Euler equation. 
\end{ilimit}

\noindent{\bf Proof.}
The ball $B$ defined in \eqref{boundsupp} endowed with the weak topology is a complete separable metrizable compact space (\cite{DS}). By \eqref{boundsupp}, we have $\supp\mu^\nu\subset B$, and thus $\mu^\nu$  satisfy the sufficient condition of Theorem \ref{Prokhorov}. Therefore there exists a subsequence $\mu^\nu$ that converges weakly in $B$ to a Borel probability measure $\mu^0$ on $B$. Because $B$ is weakly closed in $L^2(\mathbb{R}^2)$, we can extend the measure 
$\mu^0$ to $L^2(\mathbb{R}^2)$ by setting $\mu^{0}(X) = \mu^{0}(X\cap B)$ 
for any Borelian set $X$. We claim that $\mu^{0}$ is a renormalized statistical solution of the  damped and driven Euler equation. Indeed, for any $\Psi\in{\mathcal T}$, for each $i=1,2, 3$,
$$
\lim_{\nu\to 0}\int F_i(\omega)d\mu^{\nu}(\omega) = \int F_i(\omega)d\mu^{0}(\omega)
$$
holds in view of Lemma {\ref{psiregs}} because each $F_i$ is bounded and {\em weakly} continuous. In particular, the sequence $\int F_2(\omega)d\mu^{\nu}(\omega)$ is bounded, and so 
$$
\lim_{\nu\to 0}\nu \int F_2(\omega)d\mu^{\nu}(\omega) = 0
$$
holds. The fact that $\mu^{\nu}$ are SSSNS implies by Definition (\ref{sssv2}), (2), 
$$
\int \left (F_1(\omega) + F_3(\omega)\right )d\mu^{\nu}(\omega) = -\nu\int F_2(\omega)d\mu^{\nu}(\omega).
$$
Passing to the limit $\nu\to 0$ we deduce
$$
\int \left (F_1(\omega) + F_3(\omega)\right )d\mu^{0}(\omega) = 0
$$
which is the condition (\ref{renor}). Hence $\mu^0$ is a renormalized stationary statistical solution of the damped and driven Euler equation.$\quad\quad\Box$ 
\medskip

We consider the sets 
\[
B_p^\infty(r)=\left\{\omega\in B;\;\left \|\omega\right \|_{L^p({\mathbb R}^2)}\le r,\quad \left\|\omega\right\|_{L^\infty(\mathbb{R}^2)}\leq r\right\}.
\]
defined for $r>0, 1\le p< 2.$
\begin{corProk}
\label{corProk}
Let $\left\{\mu^\nu\right\}$ be a sequence of stationary statistical solutions of the damped and driven NSE in vorticity phase space, with $\nu\rightarrow 0$. Assume that  there exists $1<p<2$ and $r>0$ such that
$$
\supp\mu^\nu\subset B_p^\infty(r).
$$
Then, the limit $\mu^{0}$ of any weakly convergent subsequence 
 is a renormalized stationary statistical solution of the damped and driven Euler equation (\ref{renor}) that is supported in $B^{\infty}_p(r)$ and  satisfies the enstrophy balance (\ref{bal}).
\end{corProk}
\noindent{\bf Proof.} 
The set
$B_p^\infty (r)$ is weakly closed in $B$. Indeed, if $\omega_j\in B^{\infty}(r)$
converges weakly to $\omega$ and if $\phi\in C_0^{\infty}({\mathbb R}^2)$ then
$\left |\int\omega \phi dx\right | = \lim_{j\to\infty}\left |\int \omega_j\phi dx\right| \le r\|\phi\|_{L^1({\mathbb R}^2)}$ implies  that $\|\omega\|_{L^{\infty}({\mathbb R}^2)} \le r$.
Similarly, we obtain $\left|\int\omega \phi dx\right| \le r\|\phi\|_{L^{p'}({\mathbb R}^2)}$
where $p'>2$ is the dual exponent, and deduce that $\|\omega\|_{L^p({\mathbb R}^2)}\le r$.  
By Theorem \ref{ilimit}, the limit $\mu^{0}$ of a weakly convergent subsequence is a Borel probability measure supported in $B$ and a renormalized statistical solution of the damped and driven Euler equation (\ref{renor}). The set $U=L^{2}({\mathbb R}^2)\setminus B_p^{\infty}(r)$ is weakly open and
$\mu^{0}(U) \le \lim\inf_{\nu\to 0}\mu^{\nu}(U) = 0$ follows by general properties of weak convergence. Thus, the support of $\mu^{0}$ is included in $B_p^{\infty}(r)$.

In order to prove the enstrophy balance we consider the function 
\[
\psi^{(m)}(a_1,\ldots,a_m)=\frac{1}{2}\sum_{k=1}^{m}\left|a_k\right|^2.
\]	
Let $\left\{{\mathbf{w}}_j\right\}$ be a complete orthonormal basis in $L^2({\mathbb R}^2)$, formed with functions ${\mathbf{w}}_j\in C_0^2({\mathbb R}^2)$. Then, for each fixed $m$,
 \[
 \Psi^{(m,\epsilon)}(\omega)=\psi^{(m)}(\langle(\beta(\omega_{\epsilon}))_{\epsilon},{\mathbf{w}}_1\rangle,\ldots,\langle(\beta(\omega_{\epsilon}))_{\epsilon},{\mathbf{w}}_m\rangle)
 \]
is a function in $\mathcal{T}$, and 
\be
\langle(\Psi^{(m,\epsilon)})'(\omega),(\gamma\omega-g)\rangle=\sum_{j=1}^{m}\langle(\beta(\omega_{\epsilon}))_{\epsilon},{\mathbf{w}}_j\rangle\langle((\beta'(\omega_{\epsilon}))(\gamma\omega-g)_\epsilon)_\epsilon,{\mathbf{w}}_j\rangle.
\ee
Because $\left\{{\mathbf{w}}_j\right\}$ is an orthonormal basis in $L^2({\mathbb R}^2)$,
it follows by Parseval's theorem that
$$
\lim_{m\to\infty}\langle (\Psi^{m,\epsilon})'(\omega), (\gamma\omega-g)\rangle
= \langle (\beta(\omega_{\epsilon}))_{\epsilon}, (\beta'(\omega_{\epsilon})(\gamma\omega-g)_{\epsilon})_{\epsilon}
$$ 
holds for each $\omega$ in $B$. Moreover, because the functions 
$(\beta(\omega_{\epsilon}))_{\epsilon}$ and \\ $(\beta'(\omega_{\epsilon})(\gamma\omega-g)_{\epsilon})_{\epsilon}$ are bounded in $L^2({\mathbb R}^2)$ as $\omega\in B$, it follows that the sequence \\ $\langle (\Psi^{m,\epsilon})'(\omega), (\gamma\omega-g)\rangle$ is bounded uniformly for $\omega\in B$. Thus, we may apply the Lebesgue dominated convergence theorem to deduce
\be
\lim_{m\to\infty}\int\langle (\Psi^{m,\epsilon})'(\omega), (\gamma\omega-g)\rangle d\mu^{0}(\omega) = 
\int \langle (\beta(\omega_{\epsilon}))_{\epsilon}, (\beta'(\omega_{\epsilon})(\gamma\omega-g)_{\epsilon})_{\epsilon}d\mu^{0}(\omega).
\la{limmone}
\ee
Because $(u\cdot\nabla \omega)_{\epsilon} = \partial_k(u_k\omega)_{\epsilon}$ we have
\be
\begin{aligned}
&\langle (\Psi^{(m,\epsilon)})'(\omega),u\cdot\nabla\omega\rangle\\
&= \sum_{j=1}^{m}\langle(\beta(\omega_{\epsilon}))_{\epsilon},{\mathbf{w}}_j\rangle\langle(\beta'(\omega_{\epsilon})\partial_k(u_k\omega)_\epsilon)_\epsilon,{\mathbf{w}}_j\rangle.
\end{aligned}
\ee
In order to establish the pointwise limit
$$
\lim_{m\to\infty}\langle (\Psi^{(m,\epsilon)})'(\omega),u\cdot\nabla\omega\rangle
= \langle (\beta(\omega_{\epsilon}))_{\epsilon}, (\beta'(\omega_{\epsilon})\partial_k(u_k\omega)_{\epsilon})_{\epsilon}\rangle
$$
and the uniform bounds on $(\beta(\omega_{\epsilon}))_{\epsilon}$ and
$(\beta'(\omega_{\epsilon})\nabla_x(u\omega)_{\epsilon})_{\epsilon}$
in $L^2({\mathbb R}^2)$ for $\omega\in B$ we need to split the
Biot-Savart expression
$$
u = \frac{1}{2\pi}\frac{x^{\perp}}{|x|^2}\star \omega
$$
in two pieces, corresponding to  
$$
\ba
\frac{1}{2\pi}\frac{x^{\perp}}{|x|^2} = K_1(x) + K_2(x),\\
K_1(x) = \frac{1}{2\pi}\frac{x^{\perp}}{|x|^2}{\mathbf{1}}_{|x|\le 1},\\
K_2(x) =  \frac{1}{2\pi}\frac{x^{\perp}}{|x|^2}{\mathbf{1}}_{|x|>1}.\\
\ea
$$
Clearly, because each component of $K_1\in L^1({\mathbb R}^2)$, it follows that $u_1 = K_1\star\omega$ is in $L^2({\mathbb R}^2)$ by the Hausdorff-Young inequality, and its norm in $L^2$ is bounded by a constant
uniformly for $\omega\in B$. On the other hand, because each component of 
$K_2\in L^{p'}({\mathbb R}^2)$, with $p'>2$ the dual exponent of $p<2$, we have that $u_2 = K_2\star \omega \in L^{\infty}({\mathbb{R}^2})$ with norm bounded by $r$, as long as
$\omega \in B^{\infty}_p(r)$. Therefore $u\otimes \omega \in L^1({\mathbb R}^2)+ L^{2}({\mathbb R}^2)$, with norm bounded uniformly for $\omega\in B^{\infty}_p(r)$. Consequently,  $(u\otimes\omega)_{\epsilon}\in L^{\infty}({\mathbb R}^2)$ with norm uniformly bounded for $\omega\in B^{\infty}_p(r)$ and, because $\beta'(\omega_{\epsilon})$ is uniformly bounded in $L^2({\mathbb R}^2)$, 
we may use the Lebesgue dominated convergence
 theorem to deduce
\be
\lim_{m\to\infty}\int \langle(\Psi^{(m,\epsilon)})'(\omega),u\cdot\nabla\omega\rangle d\mu^{0}(\omega) = \int\langle (\beta(\omega_{\epsilon}))_{\epsilon}, (\beta'(\omega_{\epsilon})\partial_k(u_k\omega)_{\epsilon})_{\epsilon}\rangle
d\mu^{0}(\omega).
\la{limmtwo}
\ee
Because of (\ref{renor}, \ref{limmone}, \ref{limmtwo}) we have then
\be
\label{Repestimate}
\begin{aligned}
&\int_{L^2(\mathbb{R}^2)} \langle(\beta(\omega_{\epsilon}))_{\epsilon},(\beta'(\omega_{\epsilon}) (\gamma\omega-g)_\epsilon)_\epsilon\rangle d\mu^{0}(\omega)\\
&+\int_{L^2({\mathbb R}^2)}\langle(\beta(\omega_{\epsilon}))_{\epsilon}, (\beta'(\omega_{\epsilon})\partial_k(u_k\omega)_\epsilon)_\epsilon\rangle d\mu^0(\omega)=0.
\end{aligned}
\ee
Now we are going to investigate the term
$$
I_{\beta,\epsilon} =\int_{\mathbb R^2}(\beta(\omega_{\epsilon}))_{\epsilon}\left[\beta'(\omega_{\epsilon})\partial_k(u_k\omega)_{\epsilon}\right]_{\epsilon}dx.
$$
Integrating by parts we write
$$
I_{\beta,\epsilon} = J_{\beta, \epsilon} + K_{\beta, \epsilon},
$$
with
$$
J_{\beta, \epsilon} = -\int_{\mathbb R^2}\partial_k(\beta(\omega_{\epsilon}))_{\epsilon}
\left[\beta'(\omega_{\epsilon}) (u_k\omega)_{\epsilon}\right]_{\epsilon}dx,
$$
and
$$
K_{\beta,\epsilon} =  -\int_{\mathbb R^2}(\beta(\omega_{\epsilon})_{\epsilon})
\left[\beta''(\omega_{\epsilon})(\partial_k\omega_{\epsilon}) (u_k\omega)_{\epsilon}\right]_{\epsilon}dx.
$$
We split $J_{\beta,\epsilon}$ further, using (\ref{flux}, \ref{rho}):
$$
J_{\beta,\epsilon} = L_{\beta, \epsilon} + M_{\beta,\epsilon},
$$
with
$$
L_{\beta, \epsilon} = -\int_{\mathbb R^2}\partial_k(\beta(\omega_{\epsilon}))_{\epsilon}
\left[\beta'(\omega_{\epsilon}) (u_k)_{\epsilon}(\omega)_{\epsilon}\right]_{\epsilon}dx,
$$
and
$$
M_{\beta,\epsilon} = -\int_{\mathbb R^2}\partial_k(\beta(\omega_{\epsilon}))_{\epsilon}
\left[\beta'(\omega_{\epsilon}) \rho_{\epsilon}(u_k,\omega)\right]_{\epsilon}dx.
$$
We estimate 
$$
|M_{\beta,\epsilon}| \le C\sup|\beta|\sup|\beta'|\frac{1}{\epsilon}\|\rho_{\epsilon}(u,\omega)\|_{L^{1}({\mathbb R}^2)}.
$$
We used the fact that 
$$
\|\partial_k (\beta)_{\epsilon}\|_{L^{\infty}({\mathbb R}^2)}\le C{\frac{1}{\epsilon}}\|\beta\|_{L^{\infty}({\mathbb R}^2)}.
$$
We claim that
$$
|M_{\beta,\epsilon}| \le C\sup|\beta|\sup|\beta'|\|\omega\|_{L^2}\int_{\mathbb R^2}j(z)(1+|z|)\|\delta_{\epsilon z}\omega\|_{L^2}dz,
$$
where $(\delta_h\omega)(x) = \omega (x-h)-\omega(x)$. Indeed this follows from a bound on $\rho_{\epsilon}(u,\omega)$ and the uniform bound
$\|\delta_{\epsilon z}u\|_{L^2({\mathbb R}^2)}\le \epsilon |z| \|\omega\|_{L^2({\mathbb R}^2)}$.

We  fix $\epsilon >0$ and we consider  a sequence of compactly supported 
functions $\beta(y)$ that converge uniformly on the compact $R_{\infty} = [-2\frac{\|g\|_{L^{\infty}}}{\gamma}, 2\frac{\|g\|_{L^{\infty}}}{\gamma}]$
together with two derivatives to the function $y$ , (i.e. $\beta \to y$, $\beta'\to 1$, $\beta''\to 0$) and such that 
$$
|\beta(y) |+ |\beta'(y)| + |\beta''(y)|) \le C.
$$ 
It is easy to see that for fixed $\epsilon>0$
$$
\lim_{\beta\to y}\int_{L^2({\mathbb R}^2)} (L_{\beta, \epsilon} + K_{\beta,\epsilon})d\mu^{0}(\omega) = 0.
$$
Indeed, $K_{\beta, \epsilon}(\omega)$ is a continuous function of $\omega\in L^2({\mathbb R}^2)$, uniformly bounded on $K^{\infty}_p$ and converging pointwise to zero. As for $L_{\beta,\epsilon}$, it is also continuous, bounded and converges to $0= \int_{\mathbb R^2}\partial_k(\omega_{\epsilon})_{\epsilon} (u_k)_{\epsilon}(\omega_{\epsilon})_{\epsilon}dx$.

On the other hand, from
$$
\int_{L^2({\mathbb R}^2)} |M_{\beta,\epsilon}|d\mu^{0}(\omega) \le C\int_{L^2({\mathbb R}^2)}\int_{\mathbb R^2} j(z)(1+|z|)\|\delta_{\epsilon z}\omega\|_{L^2({\mathbb R}^2)}dz d\mu^{(0)}(\omega),
$$
with $C$ uniform for all $\beta$ in the sequence, it follows from the Lebesgue dominated convergence theorem that
$$
\lim_{\epsilon\to 0}\lim\sup_{\beta\to y}\int_{L^2({\mathbb R}^2)} |M_{\beta,\epsilon}|d\mu^{(0)} = 0.
$$
By \eqref{Repestimate} and the estimates above it follows that
$$
\lim_{\epsilon\to 0}\lim\sup_{\beta\to y}\int_{L^2(\mathbb{R}^2)} \langle(\beta(\omega_{\epsilon}))_{\epsilon},(\beta'(\omega_{\epsilon}) (\gamma\omega-g)_\epsilon)_\epsilon\rangle d\mu^0(\omega)= 0
$$
On the other hand, by the Lebesgue dominated convergence theorem again,
$$
\ba
\lim_{\epsilon\to 0}\lim\sup_{\beta\to y}\int_{L^2(\mathbb{R}^2)} \langle(\beta(\omega_{\epsilon}))_{\epsilon},(\beta'(\omega_{\epsilon}) (\gamma\omega-g)_\epsilon)_\epsilon\rangle d\mu^0(\omega) \\=
\int_{L^2(\mathbb{R}^2)}\left\{\gamma\left\|\omega\right\|_{L^2({\mathbb R}^2)}^2-\langle\mathbf{g},\omega\rangle\right\}d\mu^{0}(\omega)
\ea
$$
which proves (\ref{bal}).

\section{Long time averages and the inviscid limit}

In this section we consider SSSNSs obtained as generalized (Banach) limits of long time averages of functionals of deterministic solutions of the damped and driven Navier-Stokes equations. These SSSNS have good enough properties to pass to the inviscid limit and are used to prove that the time averaged enstrophy dissipation vanishes in the zero viscosity limit.     
\begin{Glimit}
\label{Glimit}
A generalized limit (Banach limit) is a linear continuous functional 
$$
Lim_{t\to\infty}:\mathcal{BC}([0,\infty)) \to {\mathbb R}
$$
such that 
 \begin{enumerate}
 \item
 $Lim_{t\to\infty}(g) \geq 0;\quad \forall g\in\mathcal{BC}([0,\infty))$ with $g(s)\geq0$ $\forall s\geq 0$,
 \item
 $Lim_{t\to\infty} (g) =\lim_{t\rightarrow\infty} g(t)$, 
 whenever the usual limit exists.
  \end{enumerate}
\end{Glimit}
The space $BC[0,\infty)$ is the Banach space of all bounded continuous functions defined on $[0,\infty)$, and the functional $Lim_{t\to\infty} $ is constructed as an easy application of the Hahn-Banach theorem. 
It can be shown that any generalized limit satisfies
 \be
 \label{glprop}
 \liminf_{T\rightarrow\infty}g(T)\leq Lim_{t\to\infty} (g)\leq\limsup_{T\rightarrow\infty}g(T),\quad \forall g\in\mathcal{BC}([0,\infty)).
 \ee
 Furthermore, given a particular $g_0\in\mathcal{BC}([0,\infty))$, and a sequence $t_j\rightarrow\infty$ for which $g_0(t_j)$ converges to a number $l$, we can construct a generalized limit $Lim_{t\to\infty} $ satisfying $Lim_{t\to\infty}(g_0)=l$, see \cite{bercovici, fmrt}. This implies that one can choose a functional $Lim_{t\to\infty}$ so that $Lim_{t\to\infty}g_0 = \lim\sup_{t\to\infty}g_0(t)$.

\begin{thm}\label{kb} Let $u_0\in L^2({\mathbb R}^2)$ and $\nabla^{\perp}u_0 = \omega_0\in L^{1}({\mathbb R}^2)\cap L^{\infty}({\mathbb R}^2)$. Let $f\in W^{1,1}({\mathbb {R}^2})\cap W^{1, \infty}({\mathbb {R}}^2)$. Let $t_0>0$. Let $Lim_{t\to\infty}$ be a Banach limit. Then

\be
\int_{L^2({\mathbb{R}^2})}\Phi(\omega)d\mu^{\nu}(\omega) =
Lim_{t\to\infty}\frac{1}{t}\int_{0}^{t}\Phi(S^{NS,\gamma}(s+t_0))(\omega_0))ds
 \la{limdef}
\ee
is a statistical stationary solution of the damped and driven Navier-Stokes equations. For any $p>1$ there exists $r$ depending only on $\gamma, f, \omega_0$ but not $\nu$ nor $t_0$ such that
\be
supp \mu^{\nu}\subset B_p^{\infty}(r).
\la{suppb}
\ee
The inequality
\be
\nu \int_{L^2({\mathbb R}^2)}\|\nabla \omega\|^2_{L^2({\mathbb R}^2)}d\mu^{\nu}(\omega) \le \int_{L^2({\mathbb R}^2)}\left [\langle g, \omega\rangle - \gamma\|\omega\|^2_{L^2}\right ]d\mu^{\nu}(\omega)
\la{gineq}
\ee
holds.
\end{thm}
\noindent{\bf{Proof}}. By Theorem \ref{compindi}, the set
$$
O^{+}(t_0, \{\omega_0\}) = cl\{ S^{NS,\gamma}(s+t_0)(\omega_0), \;\left | s\ge 0\,\; \right. \}
$$
is compact in $L^2({\mathbb{R}}^2)$. By Theorem {\ref{indi}}, $\Phi(S^{NS,\gamma}(s+t_0)(\omega_0))$ is a continuous bounded function on $[0,\infty)$ and so is its time average on $[0, t]$. Thus, the 
generalized limit 
$$
Lim_{t\to\infty}\frac{1}{t}\int_{0}^{t}\Phi(S^{NS,\gamma}(s+t_0))(\omega_0))ds
$$
exists. Moreover, it is a positive functional on ${\mathcal C}\left (O^{+}(t_0, \{\omega_0\})\right)$. Because of the Riesz representation theorem on compact spaces, there exists a Borel measure $\mu^{\nu}$ on the compact $O^{+}(t_0, \{\omega_0\})$ that represents the limit. The measure $\mu^{\nu}$ is supported in $O^{+}(t_0, \{\omega_0\})$,
$\mu^{\nu}(X) = \mu^{\nu}(X\cap O^{+}(t_0, \{\omega_0\}))$, for any $X$ Borelian in $L^2({\mathbb R}^2)$. 
We take a test function $\Psi\in{\mathcal T}$. Then
$$
\ba
\int_{L^2({\mathbb R}^2)}\langle \Psi'(\omega), u\cdot\nabla\omega + \gamma \omega -\nu\Delta\omega\rangle d\mu^{\nu}(\omega) \\ =
Lim_{t\to\infty} \frac{1}{t}\int_0^t\frac{d}{ds}\Psi(S^{NS,\gamma}(s+t_0)(\omega_0))ds = 0
\ea
$$
holds. This verifies definition {\ref{sssv2}} (2). In order to verify
conditions (1) and (3) we take the solution $\omega(t) =
S^{NS,\gamma}(t)(\omega_0)$ mollify it, $\omega_{\epsilon}(t) = J_{\epsilon}(\omega(t))$ and take the enstrophy balance. We obtain from (\ref{omeq})
\be
\ba
\frac{d}{2dt}\|\omega_{\epsilon}(t)\|^2_{L^2({\mathbb R}^2)} + \nu \|\nabla\omega_{\epsilon}(t)\|^2_{L^2({\mathbb R}^2)} + \gamma\|\omega_{\epsilon}(t)\|^2_{L^2({\mathbb R}^2)} - \langle g_{\epsilon}, \omega_{\epsilon}(t)\rangle  \\ = 
\langle\rho_{\epsilon}(u(t), \omega(t)), \nabla\omega_{\epsilon}(t)\rangle
\ea
\la{eebal}
\ee
Integrating in time we deduce
\be
\ba
\frac{1}{t}\int_0^t\left [\gamma\|\omega_{\epsilon}(s+t_0)\|^2_{L^2({\mathbb R}^2)} -
\langle g_{\epsilon}, \omega_{\epsilon}(s+t_0)\rangle \right ]ds +
\frac{\nu}{t}\int_0^t\|\nabla\omega_{\epsilon}(s+t_0)\|^2_{L^2({\mathbb R}^2)}ds\\
= \frac{1}{2t}\left [\|\omega_{\epsilon}(t_0)\|_{L^2({\mathbb R}^2)}^2  -  \|\omega_{\epsilon}(t+t_0)\|_{L^2({\mathbb R}^2)}^2\right ] \\
+  \frac{1}{t}\int_0^t\langle\rho_{\epsilon}(u(s+t_0), \omega(s+t_0)), \nabla\omega_{\epsilon}(s+t_0)\rangle ds
\ea 
\la{finiti}
\ee
Fixing $\epsilon>0$, we may apply $Lim_{t\to\infty}$. 
$$
\ba
Lim_{t\to\infty}\frac{1}{t}\int_0^t\left [\gamma\|\omega_{\epsilon}(s+t_0)\|^2_{L^2({\mathbb R}^2)} -
\langle g_{\epsilon}, \omega_{\epsilon}(s+t_0)\rangle \right ]ds =\\
\int_{L^2({\mathbb{R}}^2)} \left [\gamma\|\omega_{\epsilon}\|^2_{L^2({\mathbb R}^2)} -
\langle g_{\epsilon}, \omega_{\epsilon}\rangle \right]d\mu^{\nu}(\omega)
\ea
$$
and 
$$
\ba
Lim_{t\to\infty}\frac{1}{t}\int_0^t\|\nabla\omega_{\epsilon}(s+t_0)\|^2_{L^2({\mathbb R}^2)}ds = \\
\int_{L^2({\mathbb{R}}^2)}\|\nabla\omega_{\epsilon}\|^2_{L^2({\mathbb R}^2)}d\mu^{\nu}(\omega) 
\ea
$$ 
hold because the functionals are continuous. From (\ref{finiti}) we have
\be
\ba 
\int_{L^2({\mathbb{R}}^2)} \left [\gamma\|\omega_{\epsilon}\|^2_{L^2({\mathbb R}^2)} -
\langle g_{\epsilon}, \omega_{\epsilon}\rangle \right]d\mu^{\nu}(\omega) + \nu \int_{L^2({\mathbb{R}}^2)}\|\nabla\omega_{\epsilon}\|^2_{L^2({\mathbb R}^2)}d\mu^{\nu}(\omega)  \\
= Lim_{t\to\infty}\frac{1}{t}\int_0^t\langle\rho_{\epsilon}(u(s+t_0), \omega(s+t_0)), \nabla\omega_{\epsilon}(s+t_0)\rangle ds
\ea
\la{balan}
\ee
We estimate the right-hand side taking
$\nabla\omega_{\epsilon}$ in $L^{\infty}({\mathbb R}^2)$, where it costs $\epsilon^{-1}\Omega$ where $\Omega$ is a time independent bound on $\|S^{NS,\gamma}(s+t_0)\omega_0)\|_{L^{\infty}({\mathbb R}^2)}$ (from Theorem \ref{indi}). 
Then we are left with
$$
\ba
\left |Lim_{t\to\infty}\frac{1}{t}\int_0^t\langle\rho_{\epsilon}(u(s+t_0), \omega(s+t_0)), \nabla\omega_{\epsilon}(s+t_0)\rangle ds\right| \\ \le 
\Gamma Lim_{t\to\infty}\frac{1}{t}\int_0^t\int_{{\mathbb R}^2} j(z)\|\delta_{\epsilon z}\omega(s+t_0)\|_{L^2({\mathbb R}^2)}ds
\ea
$$
where $\Gamma$ is a bound on $\sup_{s\ge0}\|\omega (s+t_0)\|_{L^{\infty}({\mathbb R}^2)}\|\omega(s+t_0)\|_{L^2({\mathbb R}^2)}$.
We use crucially now the fact that $O^{+}(t_0, \{\omega_0\})$ is compact in $L^2({\mathbb{R}}^2)$. Then for every small number $h >0$ there exists $\epsilon>0$ so that
$$
\|\delta_{\epsilon z}\omega(s+t_0)\|_{L^2({\mathbb R}^2)} \le h
$$
holds {\em for all } $s\ge 0$, and all $z$ in the compact support of $j$. Therefore we have from (\ref{balan})
\be
\ba
\left |\int_{L^2({\mathbb{R}}^2)} \left [\gamma\|\omega_{\epsilon}\|^2_{L^2({\mathbb R}^2)} -
\langle g_{\epsilon}, \omega_{\epsilon}\rangle \right]d\mu^{\nu}(\omega) + \nu \int_{L^2({\mathbb{R}}^2)}\|\nabla\omega_{\epsilon}\|^2_{L^2({\mathbb R}^2)}d\mu^{\nu}(\omega)\right |  \\
\le h(\epsilon)
\ea
\la{balanc}
\ee
with $0\le h(\epsilon)$, a function satisfying $\lim_{\epsilon\to 0}h(\epsilon) = 0$. We remove now the mollifier, carefully.
First we note that
$$
\ba
\int_{L^2({\mathbb R}^2)}\left[\gamma \|\omega\|^2_{L^2({\mathbb R}^2)} - \langle g, \omega\rangle\right]d\mu^{\nu}(\omega) \\
= \lim_{\epsilon\to 0}\int_{L^2({\mathbb{R}}^2)} \left [\gamma\|\omega_{\epsilon}\|^2_{L^2({\mathbb R}^2)} - \langle g_{\epsilon}, \omega_{\epsilon}\rangle \right]d\mu^{\nu}(\omega)
\ea
$$
holds trivially because $\mu^{\nu}$ is a Borel measure. This, together with (\ref{balanc}) implies that
$$
\ba
\nu \lim\sup_{\epsilon\to 0}\int_{L^2({\mathbb{R}}^2)}\|\nabla\omega_{\epsilon}\|^2_{L^2({\mathbb R}^2)}d\mu^{\nu}(\omega) \\
\le -\int_{L^2({\mathbb R}^2)}\left[\gamma \|\omega\|^2_{L^2({\mathbb R}^2)} - \langle g, \omega\rangle\right]d\mu^{\nu}(\omega)
\ea
$$
which implies, by Fatou's lemma
\be
\ba
\nu \int_{L^2({\mathbb{R}}^2)}\|\nabla\omega\|^2_{L^2({\mathbb R}^2)}d\mu^{\nu}(\omega) \\
\le -\int_{L^2({\mathbb R}^2)}\left[\gamma \|\omega\|^2_{L^2({\mathbb R}^2)} - \langle g, \omega\rangle\right]d\mu^{\nu}(\omega).
\ea
\la{balance}
\ee
Because the right-hand side is finite, this proves (1) and (\ref{gineq}).  The proof of (3) for arbitrary $E_1, E_2$ follows from a  very similar computation as the one above. We take  $\chi'(y)$, a smooth, nonnegative,
compactly supported function defined for $y\ge 0$. 
Then $\chi(y) = \int_0^y\chi'(e)de$ is bounded on
 ${\mathbb{R}}_+$ and 
$$
\frac{d}{dt}\chi(\|\omega_{\epsilon}(t)\|^2_{L^2({\mathbb R}^2)}) =
\chi'(\|\omega_{\epsilon}(t)\|^2_{L^2({\mathbb R}^2)})\frac{d}{dt}\|\omega_{\epsilon}(t)\|^2_{L^2({\mathbb R})}.
$$ 
We multiply (\ref{eebal}) by $2\chi'(\|\omega_{\epsilon}(t)\|^2_{L^2({\mathbb{R}}^2)})$ and we proceed as above by taking time average,
long time limit and removing the mollifier. We obtain 
$$
\int_{L^2({\mathbb R}^2)}\chi'(\|\omega\|^2_{L^2({\mathbb{R}}^2)})\left\{\nu\|\nabla\omega\|^2_{L^2({\mathbb R}^2)} + \gamma \|\omega\|^2_{L^2({\mathbb R}^2)}  - \langle g, \omega\rangle\right\}d\mu^{\nu}(\omega) \le 0
$$
and letting $\chi'(y) \to {\mathbf 1}_{[E_1^2, E_2^2]}$ pointwise, with 
$0\le \chi'(y)\le 2$, concludes the proof. $\quad\Box$

\begin{thm}
\label{invlimit}
Let $f\in W^{1, 1}(\mathbb{R}^2)\cap W^{1,\infty}(\mathbb{R}^2)$. Let $u_0\in L^2(\mathbb{R}^2)$ be divergence-free and let  $\nabla^{\perp} u_0 = \omega_0\in  L^1(\mathbb{R}^2)\cap L^{\infty}({\mathbb R}^2)$. Let $\omega^{\nu}(t) = S^{NS,\gamma}(t)(\omega_0)$ be the vorticity of the solution of the damped and 
driven Navier-Stokes equation. Then,
\be
\lim_{\nu\rightarrow 0}\nu\left (\limsup_{t\rightarrow\infty}\frac{1}{t}\int_{0}^t\left\|\nabla\omega^{\nu}(s+t_0)\right\|_{L^2(\mathbb{R}^2)}^2ds\right )=0,
\la{pay}
\ee
holds for any $t_0>0$.
\end{thm}
\noindent{\bf Proof.}
We argue by contradiction and assume that the statement is false. Then, there exists a sequence $\nu_k\rightarrow 0$ and $\delta>0$, such that,
for each fixed $\nu_k$, there exists a sequence of times $t_j\to\infty $ (that may depend on $k$) such that
\be
\label{contradclaim}
\frac{\nu_k}{t_j}\int_{0}^{t_j}\left\|\nabla\omega^{\nu_k}(s+t_0)\right\|_{L^2(\mathbb{R}^2)}^2ds \ge \delta
\ee
holds for all $t_j\to\infty$. Because of the enstrophy balance
$$
\ba
\delta \le \frac{\nu_k}{t_j}\int_0^{t_j}\|\nabla\omega^{\nu_k}(s+t_0)\|^2_{L^2({\mathbb R}^2)} 
\\ =  
\frac{1}{t_j}\int_0^{t_j}\left [-\gamma \|\omega^{\nu_k}(s+t_0)\|^2_{L^2({\mathbb R}^2)} + \langle g, \omega^{\nu_k}(s+t_0)\rangle\right]ds \\
+ \frac{1}{2t_j}\left [\|\omega^{\nu_k}(t_0)\|^2_{L^2({\mathbb R}^2)} - \|\omega^{\nu_k}(t_0+ t_j)\|^2_{L^2({\mathbb R}^2)}\right]
\ea
$$
It follows that
\be
\lim\sup_{t\to\infty}\frac{1}{t}\int_0^t\left [-\gamma\|\omega^{\nu_k}(s + t_0)\|^2_{L^2({\mathbb R}^2)} + \langle g, \omega^{\nu_k}(s+t_0)\rangle \right ]ds \ge \delta.
\la{deltineq}
\ee
Because the function $-\gamma \|\omega\|^2_{L^2({\mathbb R}^2)} +\langle g, \omega\rangle$ is continuous on   $clO^{+}(t_0, \{\omega_0\})$, by the remark after Definition \ref{Glimit}, we can choose a 
generalized limit such that
\be
\begin{aligned}
&Lim_{t\to\infty}\frac{1}{t}\int_{0}^{t}\left [-\gamma\left\|\omega^{\nu_k}(s+t_0)\right\|_{L^2({\mathbb R}^2)}^2+ \langle g,\omega^{\nu_k}(s+t_0)\rangle \right] ds \\
&= \limsup_{t\rightarrow\infty}\frac{1}{t}\int_{0}^{t}\left [-\gamma\left\|\omega^{\nu_k}(s+t_0)\right\|_{L^2({\mathbb R}^2)}^2 + \langle g,\omega^{\nu_k}(s+t_0)\rangle\right ] ds.
\end{aligned}
\la{limlim}
\ee
Now, by Theorem {\ref{kb}}, this means that we have a SSSNS $\mu^{\nu_k}$ that satisfies (\ref{suppb}) and that also satisfies, in view of (\ref{deltineq}) and
(\ref{limlim})
\be
\int_{L^2({\mathbb R}^2)} \left\{-\gamma\|\omega\|^2_{L^2({\mathbb R}^2)} +
\langle g, \omega\rangle\right\}d\mu^{\nu_k}(\omega)\ge \delta>0
\la{deltimuneq}
\ee 
Passing to a weakly convergent subsequence we find with Theorem \ref{corProk}
that there exists a renormalized statistical solution of the damped and driven 
Euler equations $\mu^{0}$ that satisfies the enstrophy balance (\ref{bal}).

Because the function $\omega\mapsto \langle g, \omega\rangle$ is {\em weakly} continuous, we have
\be
\lim_{k\to\infty}\int_{L^2({\mathbb R}^2)} \langle g, \omega\rangle d\mu^{\nu_k}(\omega) =\int_{L^2({\mathbb R}^2)} \langle g, \omega\rangle d\mu^{0}(\omega) 
\la{eq}
\ee
On the other hand, by Fatou's lemma
\be
\gamma\int_{L^2({\mathbb R}^2)}\|\omega\|^2_{L^2({\mathbb R}^2)}d\mu^{0}(\omega)\le \gamma \lim\inf_{k\to\infty}\int_{L^2({\mathbb R}^2)}\|\omega\|^2_{L^2({\mathbb R}^2)}d\mu^{\nu_k}(\omega)
\la{faineq}
\ee
From (\ref{deltimuneq}) and (\ref{eq}) we have
\be
\gamma\lim\inf_{k\to\infty}\int_{L^2({\mathbb R}^2)}\|\omega\|^2_{L^2({\mathbb R}^2)}d\mu^{\nu_k}(\omega) \le \int_{L^2({\mathbb R}^2)} \langle g, \omega\rangle d\mu^{0}(\omega) -\delta
\la{deltzineq}
\ee
and from (\ref{faineq}) and (\ref{deltzineq}) we obtain
\be
\int_{L^2({\mathbb R}^2)} \left \{\gamma\|\omega\|^2_{L^2({\mathbb R}^2)} -\langle g, \omega\rangle\right\}d\mu^{0}(\omega) \le -\delta <0.
\la{inbal}
\ee
This is a contradiction because  (\ref{bal}) holds. Thus  
(\ref{pay}) holds. $\quad\quad\Box$

\begin{thm}\la{coro} Let $f\in W^{1, 1}(\mathbb{R}^2)\cap W^{1,\infty}(\mathbb{R}^2)$. Let $u_0\in L^2(\mathbb{R}^2)$ be divergence-free and let  $\nabla^{\perp} u_0 = \omega_0\in L^1(\mathbb{R}^2)\cap L^{\infty}({\mathbb R}^2)$. Let
$\mu^{\nu}$ be SSSNS associated to long time averages  given by (\ref{limdef})
that converge weakly as $\nu\to 0$ to a renormalized statistical solution $\mu^{0}$ of
the damped and driven Euler equation. Then
\be
\lim_{\nu\to 0}\int_{L^2({\mathbb R}^2)}\|\omega\|^2_{L^2({\mathbb R}^2)}d\mu^{\nu}(\omega) = \int_{L^2({\mathbb R}^2)}\|\omega\|^2_{L^2({\mathbb R}^2)}d\mu^{0}(\omega) 
\la{zl}
\ee
holds.
\end{thm}
\noindent{\bf Proof.} Indeed, by Theorem {\ref{corProk}} we know that $\mu^{0}$
satisfies (\ref{bal}). From (\ref{gineq}) and (\ref{eq}) we have
\be
\lim\sup_{\nu\to 0}\int_{L^2({\mathbb R}^2)}\gamma \|\omega\|^2_{L^2({\mathbb R}^2)}d\mu^{\nu}(\omega) \le \int_{L^2({\mathbb R}^2)} \langle g, \omega\rangle d\mu^{0}(\omega)
\la{lsup}
\ee
Using (\ref{bal}) we obtain
\be
\lim\sup_{\nu\to 0}\int_{L^2({\mathbb R}^2)}\gamma 
\|\omega\|^2_{L^2({\mathbb R}^2)}d\mu^{\nu}(\omega) \le 
\int_{L^2({\mathbb R}^2)}\gamma\|\omega\|^2_{L^2({\mathbb R}^2)}d\mu^{0}(\omega)\la{lsupp}
\ee
From (\ref{faineq}) we obtain (\ref{zl}). $\quad\quad\Box$

\vspace{1cm}
\noindent{\bf{Acknowledgment.}}\,\,\, The work of P.C. is partially supported by
 NSF-DMS grant 0504213.\,\, The work of F.R is partially supported by the Pronex in Turbulence, CNPq and FAPERJ. Brazil. grant number E-26/171.198/2003, and by CAPES Foundation. Brazil. grant number BEX4427/05-0. \\
 F.R. also wants to thank the Department of Mathematics of The University of Chicago for its hospitality.

\end{document}